\documentclass[reqno,11pt]{amsart}          
\usepackage{epsfig}                     
\usepackage{amscd}
\usepackage{amssymb}
\usepackage[all]{xy}
\usepackage{amsmath}
\usepackage{amsthm}
\newcommand{\tuborg}{\left\{\begin{array}{ll}}
\newcommand{\sluttuborg}{\end{array}\right.}
\newcommand{\calO}{\mathcal{O}}

\newcommand{\calZ}{\mathcal{Z}}

\newcommand{\calM}{\mathcal{M}}

\newcommand{\calA}{\mathcal{A}}

\newcommand{\calC}{\mathcal{C}}

\newcommand{\calH}{\mathcal{H}}

\newcommand{\calS}{\mathcal{S}}
\newcommand{\calT}{\mathcal{T}}

\newcommand{\bbH}{\mathbb{H}}
\newcommand{\bbZ}{\mathbb{Z}}

\newcommand{\bbQ}{\mathbb{Q}}
\newcommand{\bbR}{\mathbb{R}}

\newcommand{\bbA}{\mathbb{A}}

\newcommand{\ord}{{\rm ord}}

\newcommand{\spec}{{\rm Spec}}

\newcommand{\pic}{{\rm Pic}}

\newcommand{\ob}{{\rm Ob}}

\newtheorem{thm}{Theorem}[section]

\newtheorem{lemma}[thm]{Lemma}
\newtheorem{cor}[thm]{Corollary}
\newtheorem{prop}[thm]{Proposition}
\newtheorem{defi}[thm]{Definition}

\newtheorem*{remark}{Remark}
\newtheorem{remark2}[thm]{Remark}

\newtheorem{conj}[thm]{Conjecture}

\newtheorem{exercise}[thm]{Exercise}

\renewcommand{\hom}{{\rm Hom}}

\DeclareSymbolFont{AMSb}{U}{msb}{m}{n}
\DeclareMathSymbol{\N}{\mathbin}{AMSb}{"4E}
\DeclareMathSymbol{\Z}{\mathbin}{AMSb}{"5A}
\DeclareMathSymbol{\R}{\mathbin}{AMSb}{"52}
\DeclareMathSymbol{\Q}{\mathbin}{AMSb}{"51}
\DeclareMathSymbol{\I}{\mathbin}{AMSb}{"49}
\DeclareMathSymbol{\C}{\mathbin}{AMSb}{"43}

\begin{document}
 
\title{Algebraic cycles and Connes periodicity}
\author{Jinhyun Park}
\address{Department of Mathematics, The University of Chicago, Chicago, Illinois 60637, USA}
\email{jinhyun@math.uchicago.edu}
\date{February 14, 2006. Last modified on July 11, 2006.}
\begin{abstract}We apply the classical technique on cyclic objects of Alain Connes to various objects, in particular to the higher Chow complex of S. Bloch to prove a Connes periodicity long exact sequence involving motivic cohomology groups. The Cyclic higher Chow groups and the Connes higher Chow groups of a variety are defined in the process and various properties of them are deduced from the known properties of the higher Chow groups. Applications include an equivalent reformulation of the Beilinson-Soul\'e vanishing conjecture for the motivic cohomology groups of a smooth variety $X$ and a reformulation of the conjecture of Soul\'e on the order of vanishing of the zeta function of an arithmetic variety.
\end{abstract}

\maketitle

\tableofcontents

\vskip1cm

\section*{\textbf{Introduction}}

Let $R$ be a commutative ring with unity and $A$ an associative unitary $R$-algebra. One classical result proven by A. Connes is that the Hochschild homology groups $HH_*(A)$ and the cyclic homology groups $HC_*(A)$ of $A$ over $R$ have the Connes periodicity exact sequence
\begin{equation}
\cdots \overset{S}{\to} HC_{n-1} (A) \overset{B}{\to} HH_{n} (A) \overset{I}{\to} HC_n (A) \overset{S}{\to} HC_{n-2} (A) \overset{B}{\to} \cdots.
\end{equation} One standard way of proving this sequence is to construct a $2$-periodic cyclic bicomplex from the Hochschild complex $(A^{\otimes *+1}, b)$ and the acyclic augmented bar complex $(A^{\otimes *+1}, b')$. Furthermore, A. Connes observed that a similar construction works for any \emph{cyclic modules} and formalized it. It will be reviewed in \S1 and \S 2.

An interesting observation we will see in this paper is that various familiar theories, e.g. singular homology complex or the higher Chow complex of S. Bloch, also give cyclic modules, thus in particular we can prove an analogue of Connes periodicity exact sequence involving the higher Chow groups $CH_q (X, n)$ and their cyclic cousins, the Cyclic higher Chow groups $CHC_q (X, n)$:
\begin{eqnarray*}\cdots  CH_q (X, n)  \overset{I}{\to} CHC_q (X, n) \overset{S}{\to} CHC_{q-2} (X, n-2)  \\ \overset{B}{\to}  CH_{q-1} (X, n-1)  \overset{I}{\to} \cdots.\end{eqnarray*}

The higher Chow groups play the role of the motivic cohomology groups for smooth varieties (See also  \cite{V}) thus we in fact obtain a periodicity sequence involvoing motivic cohomology groups. This periodicity and the functoriality result show that the study of the higher Chow groups is equivalent to the study of the Cyclic higher Chow groups in the following sense: \emph{a morphism $f: X \to Y$ of varieties induces isomorphisms of the higher Chow groups if and only if it induces isomorphisms of the Cyclic higher Chow groups.} 

Further adapting the idea of A. Connes, it is possible to define the \emph{Connes higher Chow groups} which are isomorphic to the Cyclic higher Chow groups after tensoring with $\bbQ$. Thus, we can rephrase the above equivalence in the following way: \emph{a morphism $f: X \to Y$ of varieties induces isomorphisms of the rational higher Chow groups if and only if it induces isomorphisms of the rational Connes higher Chow groups.} 

S. Bloch proved a list of properties, for example, $\bbA^1$-homotopy, localization sequence, etc., for the higher Chow groups. The above equivalence then suggests that similar series of results may have their cyclic counterparts, and many of them can be proven. It also supplies a new way of looking at various conjectures on the motivic cohomology groups, for example, the Beilinson-Soul\'e vanishing conjecture. They are discussed after the proof of the Connes periodicity. 

The proof of the periodicity follows the one presented in \cite{L} and the involved techniques are no more than some basic homological algebra.

\section{Precyclic objects with the last degeneracy}

Let $\calC$ be a category. All definitions in this section except for the last one are more or less standard and taken from \cite{L}.

\begin{defi}\emph{A presimplicial object} in $\calC$ is a collection $\{C_n\}$, $n \geq 0$ of objects in $\calC$ with morphisms
\begin{equation}
d_i : C_n \to C_{n-1}, \ \ 0 \leq i \leq n
\end{equation} satisfying the property
\begin{equation}
d_i d_j = d_{j-1} d_i, \ \ i<j.
\end{equation}

\emph{A simplicial object} in $\calC$ is a presimplicial object $\{C_n \}$ in $\calC$ with more morphisms 
\begin{equation}
s_i : C_n \to C_{n+1}, \ \ 0 \leq i \leq n
\end{equation}satisfying
\begin{eqnarray}
& s_i s_j = s_{j+1} s_i, \ \ i \leq j \\
& d_i s_j = \tuborg s_{j-1} d_i, & i<j \\ \mbox{id}, & i=j, i = j+1 \\ s_j d_{i-1} , & i>j+1 \sluttuborg
\end{eqnarray}

\end{defi}

\begin{defi}\emph{A precyclic object} in $\calC$ is a presimplicial object $\{ C_n \}$ with the cyclic morphisms
\begin{equation}
T_n: C_n \to C_n
\end{equation} satisfying
\begin{eqnarray}
& T_n ^{n+1} = \mbox{id}, \\
& \tuborg d_0 T_n = d_n, & \\
 d_i T_n = T_{n-1} d_{i-1}, &1 \leq i \leq n.\sluttuborg
\end{eqnarray}

\emph{A cyclic object} in $\calC$ is a simplicial object $\{ C_n \}$ in $\calC$ that is also a precyclic object such that
\begin{eqnarray}
& s_0 T_n = T_{n+1} ^2 s_n, \\
& s_i T_n = T_{n+1} s_{i-1} , \ \  1 \leq i \leq n.
\end{eqnarray}

\end{defi}

The above concepts can be formalized in terms of functors: let $\Delta$ be the simplicial category, i.e. a category whose objects are $[n]:=\{ 0, 1, \cdots, n \}$ and morphisms are generated by faces $\delta_i : [n-1] \to [n]$, $0 \leq i \leq n$ and degeneracies $\sigma_j : [n+1] \to [n]$, $0 \leq j \leq n$, with the following properties:
\begin{eqnarray}
&\delta_j \delta_i = \delta_i \delta_{j-1}, \ \ i < j \\
& \sigma_j \sigma_i = \sigma_i \sigma_{j+1}, \ \ i \leq j \\
& \sigma_j \delta_i = \tuborg \delta_i \sigma_{j-1}, & i<j \\ \mbox{id}, & i=j, i=j+1 \\ \delta_{i-1} \sigma_j. & i>j+1 \sluttuborg
\end{eqnarray}

The cyclic category $\Delta C$ is a category with $\ob (\Delta C) = \ob (\Delta)$ and morphisms are generated by the faces, degeneracies and the cyclic operators $\tau_n : [n] \to [n]$ with the extra properties
\begin{eqnarray}
& \tau_n \delta_i = \delta_{i-1} \tau_{n-1}, \ \ 1 \leq i \leq n \\
& \tuborg \tau_n \sigma_0 = \sigma_n \tau_{n+1} ^2 & \\ \tau_n \sigma_i = \sigma_{i-1} \tau_{n+1}, & 1 \leq i \leq n \sluttuborg \\ 
& \tau_n ^{n+1} = \mbox{id}.
\end{eqnarray}

Then a simplicial object in $\calC$ is a functor
\begin{equation}
\Delta^{op} \to \calC
\end{equation} and a cyclic object in $\calC$ is a functor
\begin{equation}
\Delta C^{op} \to \calC.
\end{equation} To define the presimplicial objects and precyclic objects in terms of functors, we simply create categories $\Delta_{pre}$ and $\Delta C_{pre}$ without degeneracies.

Let's define one more terminology that encodes the minimal list of properties we need.
\begin{defi}A precyclic object in $\calC$ is said to be equipped \emph{with the last degeneracy} if there is a morphism $s_n: C_n \to C_{n+1}$ such that 
\begin{equation}
d_i s_n = \tuborg s_n d_{i-1}, & 0 \leq i \leq n-1 \\ \mbox{id}, & i = n, n+1. \sluttuborg
\end{equation}
\end{defi}

\section{Formalism of Connes periodicity}

\subsection{The cyclic bicomplex}Let $\calA$ be an abelian category. From now, $\{C_n\}$ will be a fixed precyclic object with the last degeneracy in $\calA$. Define the signed cyclic operator \begin{equation}t=t_n: C_n \to C_n\end{equation} to be $t_n = (-1)^n T_n$, where $-1$ is the additive inverse of the element $1=\mbox{id}$ in the unitary associative not necessarily commutative ring $\hom_{\calA} (C_n, C_n)$, and $(-1)^n$ is the $n$-time self composition of the endomorphism $-1$.

Then, the following is the comprehensive list of all properties of $\{C_n \}$:
\begin{enumerate}

\item $d_i d_j = d_{j-i} d_i, \ \ i<j$
\item $d_i s_n = \tuborg s_{n-1} d_{i-1}, & i<n \\ \mbox{id}, & i = n, n+1 \sluttuborg$
\item $\tuborg d_0 t_n = (-1)^n d_n & \\ d_i t_n = - t_{n-1} d_{i-1} , & 1 \leq i \leq n \sluttuborg $
\item $t_n ^{n+1} = \mbox{id}.$
\end{enumerate}

We will prove that this list is necessary and sufficient for the proof of the Connes periodicity long exact sequence. Define a map
\begin{equation}
\partial:= \sum_{i=0} ^n (-1)^i d_i: C_n \to C_{n-1}.
\end{equation} The property 1. shows that $\partial^2 = 0$ thus we have a complex
\begin{equation}
\cdots \overset{\partial}{\to} C_{n+1} \overset{\partial}{\to} C_n \overset{\partial}{\to} C_{n-1} \overset{\partial}{\to} \cdots.
\end{equation}

\begin{defi} The $n$-th homology $H_n (C_\bullet)$ of the presimplicial object $\{C_n\}$ is the homology at $C_n$ of the above complex.
\end{defi}

Define the \emph{extra degeneracy morphism} (\emph{c.f.} 1.1.12 and 2.5.7 in \cite{L}) \begin{equation}s= (-1)^{n+1} t_{n+1} s_n: C_n \to C_{n+1}\end{equation} and define a variation of the boundary operator $\partial$:
\begin{equation}
\partial ' := \sum_{i=0} ^{n-1} (-1)^i d_i : C_n \to C_{n-1}.
\end{equation} Notice that we did not add $(-1)^n d_n$. That $(\partial')^2 =0$ follows again from the property 1. An interesting essential point is to notice the following amusing interaction of $\partial'$ with the extra degeneracy $s$:

\begin{lemma}\label{acyclicity}The complex $(C_\bullet, \partial')$ is acyclic. (\emph{c.f.} 1.1.11 and 1.1.12 in \emph{ibid.})
\end{lemma}
\begin{proof} The point is to notice that the extra degeneracy gives an contracting homotopy of this complex. First, note that
\begin{equation}
d_0 s = (-1)^{n+1} d_0 t_{n+1} s_n = d_n s_n = \mbox{id},
\end{equation} where the properties 3. and 2. were used respectively. For $1 \leq i \leq n$, we have
\begin{eqnarray}
d_i s &=& (-1)^{n+1} d_i t_{n+1} s_n = (-1)^n t_n d_{i-1} s_n\\
& = & (-1)^n t_n s_{n-1} d_{i-1} = s d_{i-1}.
\end{eqnarray}

Thus,

\begin{eqnarray}
\partial' s & = & \sum_{i=0} ^n (-1)^i d_i s = \mbox{ id } + \sum_{i=1} ^n (-1)^i s d_{i-1} \\
&=& \mbox{id} - s \sum_{i=0} ^{n-1} (-1)^i d_{i} = \mbox{id} - s \partial'.
\end{eqnarray} This shows the desired assertion.
\end{proof}

\begin{remark}{\rm In fact, S. Bloch remarked in an private communication with the author that we can instead use just the last degeneracy only to show the assertion of the lemma: let $s'= (-1)^n s_{n}: C_n \to C_{n+1} $

Then, 
\begin{eqnarray}
\partial' s' &=& (-1)^n \partial'  s_n =(-1)^n \sum_{i=0} ^n (-1)^i d_i s_n \\ &=&(-1)^n \sum_{i=0} ^{n-1} (-1)^i s_{n-1} \partial_i + \mbox{id} \\
&=& - (-1)^{n-1} s_{n-1} \sum_{i=0} ^{n-1} (-1)^i d_i + \mbox{ id} = - s' \partial' + \mbox{id},
\end{eqnarray}i.e. $\partial' s' + s' \partial' = \mbox{id}$, which reproves the assertion.}\end{remark}

Now, we prove that the property 3:
\begin{equation}\label{first identity}
d_0 t_n = (-1)^n d_n
\end{equation}
\begin{equation}\label{second identity}
d_i t_n = - t_{n-1} d_{i-1}, \ \ 1 \leq i \leq n.
\end{equation}
implies that
\begin{equation}
\tuborg \partial (1-t) = (1-t)\partial ' \\ \partial ' N = N \partial \sluttuborg
\end{equation}
where $N= N_n: C_n \to C_n$ is defined as $N = 1 + t + \cdots + t^n$. Notice that the property 4. shows that $N (1-t) = (1-t) N = 0$ so that these two equations will give a bicomplex (\emph{c.f.} 2.5.5 in \emph{ibid.}) $CC(C_\bullet)$ that sits in the first quadrant and in particular the rows stop on the left in degree $0$ and the bottom left corner is in bidegree $(0,0)$:

$$\xymatrix{\vdots \ar[d]_{\partial} &\vdots \ar[d] _{- \partial'} & \vdots \ar[d] _{\partial} & \vdots \ar[d] _{- \partial '} & \\
C_{n+1} \ar[d] _{\partial} &  C_{n+1} \ar[l]_{1-t} \ar[d] _{- \partial'} & C_{n+1} \ar[l] _N \ar[d] _{\partial} & C_{n+1} \ar[l]_{1-t} \ar[d] _{- \partial '} & \cdots \ar[l] _N \\
C_n \ar[d] _{\partial} & C_n \ar[d] _{- \partial '} \ar[l]_{1-t} & C_n \ar[d] _{\partial} \ar[l] _{N} &  C_n \ar[l] _{1-t} \ar[d] _{- \partial'} & \cdots \ar[l]_N \\
\vdots  \ar[d] _{\partial} & \vdots \ar[d]_{- \partial '} & \vdots \ar[d] _{\partial}& \vdots  \ar[d] _{-\partial '}& \vdots \\ C_0 & C_0\ar[l] _{1-t} & C_0 \ar[l]_N & C_0 \ar[l] _{1-t} & \cdots \ar[l] _N}$$

\begin{cor}\label{first commutativity} $\partial (1-t) = (1-t) \partial '$.(\emph{c.f.} Lemma 2.1.1 in \cite{L})
\end{cor}
\begin{proof}
Note that the first $t =t_n$ and the second $t =t_{n-1}$. When $i=0$, by \eqref{first identity}, $d_0 - d_0 t = d_0 - (-1)^n d_n$, and when $1 \leq i \leq n$, by \eqref{second identity}, $d_i - d_i t = d_i - d_i t_n = d_i + t_{n-1} d_{i-1}$. Thus,
\begin{eqnarray}
\partial(1-t) & = & \sum_{i=0} ^n (-1)^i d_i (1-t)\\ & = &d_0 - (-1)^n d _n + \sum_{i=1} ^n (-1)^i \left( d_i + t_{n-1} d_{i-1} \right) \\
& = & d_0 - (-1)^n d _n + \sum_{i=1} ^n (-1)^i d _i + t_{n-1} \sum_{i=1} ^n (-1)^i d_{i-1} \\
& = & d_0 + \sum_{i=1} ^{n-1} (-1)^i d_i + t_{n-1} \sum_{i=0} ^{n-1} (-1)^{i+1} d_i \\
& = & \partial ' - t_{n-1} \partial' = (1-t)\partial '
\end{eqnarray} as desired.\end{proof}

\begin{lemma}\label{presecond commutativity} We have (\emph{c.f.} Lemma 2.1.1 in \emph{ibid.})
\begin{equation}\label{third identity}
d_i t^j = (-1)^j t^j d_{i-j}, \ \ \mbox{ for } i \geq j,
\end{equation}
\begin{equation}\label{fourth identity}
d_i t^j = (-1)^{n-j +1} t^{j-1} d_{n+1+i-j}, \ \ \mbox{ for } i < j.
\end{equation}
\end{lemma}
\begin{proof}Notice that the first $t=t_n$ and the second $t= t_{n-1}$. When $i \geq j$, by applying \eqref{first identity} successively, we obtain
\begin{eqnarray}
d_i  t^j & = & d_i t_n ^j = - t_{n-1} d_{i-1} t_n ^{j-1} \\
& = & (-1)^2 t_{n-1} ^2 d_{i-2} t_n ^{j-2} \\
& = & \cdots = (-1)^j t_{n-1} ^j d_{i-j} = (-1)^j t^j d_{i-j},
\end{eqnarray} thus obtain \eqref{third identity}.

When $i<j$, by appying \eqref{second identity} successively,
\begin{eqnarray}
d_i t^j & = & d _i t_n ^j = - t_{n-1} d_{i-1} t_n ^{j-1} \\
& = & (-1)^2 t_{n-1} ^2 d_{i-2} t_n ^{j-2}\\
& = & \cdots = (-1)^i t_{n-1} ^i d_0 t_n ^{j-i}\\
& \overset{\mbox{\eqref{first identity}}}{=} & (-1)^i t_{n-1} ^i (-1)^n d_n t_n ^{j-i-1}.
\end{eqnarray} Since $n \geq j-i-1$, by using \eqref{third identity} which was just proven above, the last expression is
\begin{eqnarray}& = & (-1)^i t_{n-1} ^i (-1)^n \cdot (-1)^{j-i-1} t_{n-1} ^{j-i-1} d_{n+1+i-j} \\
& = & (-1)^{n+j-1} t^{j-1} d _{n+1+i-j} \\
& = & (-1)^{n - j + 1} t^{j-1} d_{n+1+i-j},
\end{eqnarray} which establishes \eqref{fourth identity}.
\end{proof}

\begin{cor}\label{second commutativity}$\partial' N = N \partial$.(\emph{c.f.} Lemma 2.1.1 in \cite{L})
\end{cor}

\begin{proof}Notice that the first $N=N_n$ and the second $N=N_{n-1}$. 
\begin{eqnarray}
\partial' N &=& \left( \sum_{i=0} ^{n-1} (-1)^i d_i \right) \left( \sum_{j=0} ^n t^j \right) = \sum_{i=0} ^{n-1} \sum_{j=0} ^n (-1)^i d_i t^j\\
& = & \sum_{i \geq j} (-1)^i d_i t^j + \sum_{i<j} (-1)^i d_i t^j \\
& = & \sum_{0 \leq j \leq i \leq n-1} (-1)^{i-j} t^j d_{i-j} \\ &  & + \sum_{0 \leq i < j \leq n } (-1)^{n+1+i-j} t^{j-1} d_{n+1+i-j}.
\end{eqnarray} Let $0 \leq s \leq n$. Then the coefficient of $(-1)^s d_s$ is
\begin{equation}
\sum_{0 \leq j \leq n-1 -s} t^j + \sum_{n-s \leq j-1 \leq n-1} t^{j-1} = \sum_{j=0} ^{n-1} t^j = N_{n-1}.
\end{equation} Hence, indeed $\partial' N = N \partial.$
\end{proof}

\begin{defi}The \emph{cyclic homology} (\emph{c.f.} Definitions 2.1.3 and 2.5.6 in \emph{ibid.}) $HC_n (C_\bullet)$ of the precyclic object $\{C_n \}$ is defined to be the $n$-th homology of the total complex of the above bicomplex $CC(C_{\bullet})$
\begin{equation}
HC_n (C_{\bullet}) = H_n \left( Tot (CC(C_\bullet))\right).
\end{equation}
\end{defi}

Let $C_n ^{\lambda} = C_n /( 1-t)$, the cokernel of the endomorphism $(1-t)$. Since $\partial (1-t) = (1-t) \partial'$, $\partial$ in fact descends onto $C_n ^{\lambda}$, thus giving the \emph{Connes complex} (\emph{c.f.} 2.1.4 and 2.5.9 in \emph{ibid.}) $(C_{\bullet} ^{\lambda}, \partial)$ of the precyclic object $\{ C_n \}$:
\begin{equation}
\cdots \overset{\partial}{\to} C_{n+1} ^{\lambda} \overset{\partial}{\to} C_n ^{\lambda} \overset{\partial}{\to} C_{n-1} ^{\lambda} \overset{\partial}{\to} \cdots
\end{equation} and its homoogy
\begin{equation}
HC_n ^{\lambda} (C_{\bullet}):= H_n (C_{\bullet} ^{\lambda}, \partial)
\end{equation} is the \emph{Connes homology} of $\{ C_n \}$. From the bicomplex, we note that there is a natural morphism
\begin{equation}
HC_n (C_{\bullet}) \to HC_n ^{\lambda} (C_{\bullet}),
\end{equation} which is not always an isomorphism, but \emph{almost} an isomorphism in the following sense:

\begin{thm}[\emph{c.f.} Theorem 2.1.5 in \emph{ibid.}]The above natural morphism is an isomorphism if all nonzero integers $n$ can be inverted in $\calA$. Even though this is not possible, the induced map
\begin{equation}
HC_n (C_{\bullet}) \otimes_{\bbZ} \bbQ \to HC_n ^{\lambda} (C_{\bullet}) \otimes_{\bbZ} \bbQ
\end{equation} is an isomorphism, where the tensor product means that we invert integers formally.
\end{thm}
\begin{proof}Consider the $n$-th row of the bicomplex $CC$.  We will construct a contracting homotopy of this complex. Let $h':= \frac{1}{n+1} \mbox{id}$, $h:= - \frac{1}{n+1} \sum_{i=1} ^n i t^i$.
$$
\xymatrix{ C_n \ar[d] ^{\mbox{id}} \ar[dr]^h & C_n \ar[l] _{1-t} \ar[d]^{\mbox{id}} \ar[dr] ^{h'} & C_n \ar[d] ^{\mbox{id}} \ar[l] _{N} \ar[dr] ^{h} & C_n \ar[l] _{1-t} \ar[d]^{\mbox{id}} \ar[dr]^{h'} & \cdots \ar[l]_{N} \ar[d]^{\mbox{id}} \\
C_n & C_n \ar[l]_{1-t} & C_n \ar[l] _{N} & C_n\ar[l] _{1-t} & \cdots \ar[l] _{N}  }
$$
Then, 
\begin{eqnarray}
(n+1) \left( h'N + (1-t)h \right) & = & \mbox{id} \circ N - (1-t) \sum_{i=1} ^n i t^i \\ &=& \mbox{id} + \sum_{i=1} ^n t^i - \sum_{i=1} ^n i t^i + \sum_{i=1} ^n i t^{i+1} \\
& = & \mbox{id} + n \cdot \mbox{id} = (n+1) \mbox{id}
\end{eqnarray} thus, $h'N + (1-t) h = \mbox{id}$. Similarly we see that $Nh' + h(1-t) = \mbox{id}$. Thus, if integers are invertible, each row is an acyclic augmented complex with the $0$-th homology isomorphic to $HC_n ^{\lambda} (C_{\bullet})$. In case it is not, we invert integers first and then we have the same property. Then, the spectral sequence associated to this bicomplex shows the assertion.\end{proof}

\subsection{Connes periodicity}In this subsection we prove the Connes periodicity exact sequence for a precyclic object $\{C_n \}$ with the last degeneracy in an abelian category $\calA$. Actually this is a very simple corollary of the [2,0]-periodicity of the bicomplex $CC$. Notice first that by the Lemma \ref{acyclicity}, each odd numbered column of $CC$ is acyclic. Let $CC^{\{2 \}}$ be the bicomplex consisting of only column $0$ and $1$ of $CC$. By the Lemma \ref{acyclicity}, the natural inclusion of the complex $(C_\bullet, \partial)$ into $CC^{\{2 \}}$ as the 0-th column is a quasi-isomorphism. Now we prove the main theorem:

\begin{thm}[Periodicity (\emph{c.f.} 2.2.1 and Theorem 2.5.8 in \emph{ibid.})]Let $\{C_n \}$ be a precyclic object with the last degeneracy in an abelian category. Then there is a natural long exact sequence
\begin{equation}
\cdots \overset{B}{\to}  H_n \overset{I}{\to} HC_n \overset{S}{\to} HC_{n-2} \overset{B}{\to}  H_{n-1} \overset{I}{\to} \cdots.
\end{equation}
\end{thm}
\begin{proof}We have a short exact sequence of bicomplexes
\begin{equation}
0 \to CC ^{\{ 2 \}} \to CC \to CC [2,0] \to 0,
\end{equation}where [2,0] indicates that the degrees are shifted horizontally by $2$ and vertically by $0$. Then by taking $Tot$ and the homology long exact sequence, we obtain the periodicity sequence.
\end{proof}

\section{Examples}
In this section, we will see various examples of cyclic objects. The first one is well-known but the rest three haven't been studied from this perspective so far until now. All examples are indeed well-known to be simplicial objects so that what matters here is to check only the conditions:

\begin{eqnarray}\label{cyclic compatibility}
\tuborg d_0 T_n = d_n & \\ d_i T_n = T_{n-1} d_{i-1} & 1 \leq i \leq n \sluttuborg
\end{eqnarray} for the unsigned cyclic operator $T_n$.

\subsection{Hochschild complex}Let $k$ be a commutative unitary ring and let $A$ be a unitary $k$-algebra. Let $C_n (A):= A ^{\otimes (n+1)}$ with the face maps
\begin{eqnarray}
\tuborg d_i (a_0 \otimes \cdots \otimes a_n) = a_0 \otimes \cdots \otimes a_i a_{i+1} \otimes \cdots \otimes a_n & 0 \leq i \leq n-1 \\ d_n (a_0 \otimes \cdots \otimes a_n) = a_n a_0 \otimes a_1 \otimes \cdots \otimes a_{n-1} & i = n \sluttuborg
\end{eqnarray} and the degeneracies
\begin{equation}
s_j (a_0 \otimes \cdots \otimes a_n) = a_0 \otimes \cdots \otimes a_j \otimes 1 \otimes a_{j+1} \otimes \cdots \otimes a_n, \ \ 0 \leq j \leq n.
\end{equation} It is a simplicial $k$-module whose homology groups are, by definition, the Hochschild homology (\emph{c.f.} 1.1.3 in \emph{ibid.}) groups $HH_n (A)$. Define the unsigned cyclic operator 

\begin{equation}
T_n: A^{\otimes (n+1)} \to A^{\otimes (n+1)}
\end{equation} by $T_n (a_0 \otimes \cdots \otimes a_n) = a_n \otimes a_0 \otimes \cdots \otimes a_{n-1}$.

\begin{prop}It satisfies the cyclic compatibility conditions \eqref{cyclic compatibility}.
\end{prop}

\begin{proof}
\begin{eqnarray}
d_0 T_n (a_0 \otimes \cdots \otimes a_n) &=& d_0 (a_n \otimes a_0 \cdots \otimes a_{n-1} ) \\
&=& a_n a_0 \otimes a_1 \otimes \cdots \otimes a_{n-1} \\
&=& d_n (a_0 \otimes \cdots \otimes a_n).
\end{eqnarray} Thus, $d_0 T_n = d_n$. For $1 \leq i \leq n$,

\begin{eqnarray}
d_i T_n (a_0 \otimes \cdots \otimes a_n) &= &d_i (a_n \otimes a_0 \otimes \cdots \otimes a_{n-1} )\\
&=& a_n \otimes a_0 \otimes \cdots \otimes a_{i-1} a_i \otimes \cdots \otimes a_{n-1}\\
T_n d_{i-1} (a_0 \otimes \cdots \otimes a_n) &=& T_n (a_0 \otimes \cdots \otimes a_{i-1} a_i \otimes \cdots \otimes a_n) \\
&=& a_n \otimes a_0 \otimes \cdots \otimes a_{i-1} a_i \otimes \cdots \otimes a_{n-1},
\end{eqnarray} proving $d_i T_n = T_n d_{i-1}$.\end{proof}

Thus, by the above machine, we deduce the usual Connes peridicity involvong the usual Hochschild homology the homology of $\{C_n(A)\}$, and the cyclic homology.

\subsection{Singular complex}
It is interesting to note that this machine works for singular homology groups of topological spaces. Let $X$ be a topological space and consider the geometric simplicies
\begin{equation}
\Delta_n = \left\{ (t_0, \cdots, t_n) | 0 \leq t_i \leq 1 , \sum t_i = 1 \right\} \subset \bbR^{n+1}.
\end{equation} Define the face maps and the degeneracy maps
\begin{eqnarray}
& \delta_i : \Delta_{n-1} \to \Delta_n \\ & (t_0, \cdots, t_{n-1}) \mapsto (t_0, \cdots, t_{i-1} , 0, t_i, \cdots, t_{n-1} ); \\
& \sigma_i : \Delta_{n+1} \to \Delta_n \\  & (t_0, \cdots, t_{n+1}) \mapsto (t_0, \cdots, t_{i-1}, t_i + t_{i+1} , t_{i+2}, \cdots, t_{n+1}).
\end{eqnarray}

A singular $n$-simplex is a continuous function $f: \Delta_n \to X$. Let $S_n (X)$ be the free abelian group on the set of all singular $n$-simplices of $X$. Then the above faces and degeneracies induce the structure of a simplicial abelian group on $\{ S_n (X) \}$ as follows:

\begin{eqnarray}
d_i: S_n (X) \to S_{n-1} (X); & d_i f (u):= f(\delta_i u), \ u \in \Delta_{n-1} \\
s_i: S_n (X) \to S_{n+1} (X); & s_i f (v) := f(\sigma_i v), \ v \in \Delta_{n+1}.
\end{eqnarray}

Furthermore, we now see that this simplicial abelian group has a natural cyclic structure. (Though, we will avoid using the name ``cyclic abelian group'' for an apparent reason here.) The geometric simplex $\Delta_n$ has a natural cyclic operator $T_n: \Delta_n \to \Delta_n$ defined by \begin{equation} T_n (u_0, \cdots, u_n) := (u_n, u_0, \cdots, u_{n-1})\end{equation} for $(u_0, \cdots, u_{n}) \in \Delta_n$. This induces an unsigned cyclic operator
\begin{equation}
T_n ^* : S_n (X) \to S_n (X)
\end{equation} defined by sending a singular $n$-simplex $f \mapsto f\circ T^{-1}$.

\begin{prop}
It satisfies the cyclic compatibility conditions \eqref{cyclic compatibility}.
\end{prop}

\begin{proof}
For an arbitrary singular $n$-simplex $f$ and a point $(u_0, \cdots, u_{n-1}) \in \Delta_{n-1}$, we have
\begin{eqnarray*}
(d_0 T^* f)(u_0, \cdots, u_{n-1}) &=& f(T^{-1} \delta_0 (u_0, \cdots, u_{n-1}))\\
&=& f(T^{-1} (0, u_0, \cdots, u_{n-1})) \\
&=& f(u_0, \cdots, u_{n-1}, 0), \\
(d_n f)(u_0, \cdots, u_{n-1}) &=& f(\delta_n (u_0, \cdots, u_{n-1})) = f(u_0, \cdots, u_{n-1}, 0),
\end{eqnarray*} thus, $d_0 T^* = d_n.$ For $1 \leq i \leq n$,
\begin{eqnarray*}
(d_i T^* f)(u_0, \cdots, u_{n-1}) &=& f(T^{-1} \delta_i (u_0, \cdots, u_{n-1})) \\
&=& f(T^{-1} (u_0, \cdots, u_{i-1} , 0, u_i, \cdots, u_{n-1})) \\
&=& f(u_1, \cdots, u_{i-1} , 0, u_i, \cdots, u_{n-1}, u_0 ),\\
(T^* d_{i-1} f)(u_0, \cdots, u_{n-1}) &=& f(\delta_{i-1} T^{-1} (u_0, \cdots, u_{n-1})) \\
&=& f(\delta_{i-1} (u_1, \cdots, u_{n-1}, u_0)) \\
&=& f(u_1, \cdots, u_{i-1}, 0, u_i, \cdots, u_{n-1}, u_0),
\end{eqnarray*} thus, $d_i T^* = T^* d_{i-1}$. This proves the proposition.
\end{proof}

The homology of $S_n (X)$ is called the singular homology groups and usually denoted by $H_n (X)$. The cyclic homology of $S_n (X)$ will be denoted by $HC_n (X)$ and let's call it the cyclic singular homology groups. It seems that it is possible to prove that the cyclic singular homology groups of $X$ have various analogous properties, for example, homotopy invariance, Mayer-Vietoris sequence, to name a few. However, it has a different homology groups when $X$ is a point, for example. It is an interesting exercise to prove the following simple result using the corresponding Connes periodicity sequence:

\begin{exercise}
The cyclic singular homology groups of a point is
\begin{equation}
HC_n (\{ * \}) = \tuborg \bbZ & n: \mbox{even} \\ 0 & n: \mbox{odd}.\sluttuborg
\end{equation}
\end{exercise}

It could be an interesting result if one can relate this group to some groups already known to us well.

\subsection{Higher Chow complex}
Let's recall the definition of higher Chow groups from \cite{B1}. Let $k$ be a field and $X$ be a scheme of finite type over $k$. For $n \geq 0$, let $\Delta^n = \spec \left( k[t_0, \cdots, t_n]/ \left(\sum_{i=0} ^n t_i -1 \right) \right)$ with the faces ($0 \leq j \leq n$)

\begin{eqnarray}
\partial_j: \Delta^{n-1} \to \Delta^n ; \ \ \ \partial_j ^* (t_i) = \tuborg t_i & i<j, \\ 0 & i=j,\\ t_{i-1} & i>j.\sluttuborg
\end{eqnarray}
and the degeneracies
\begin{eqnarray}
\pi_j: \Delta^n \to \Delta^{n-1} ;\ \ \ \pi_j ^* (t_i) = \tuborg t_i & i<j, \\ t_i + t_{i+1} & i=j, \\ t_{i+1} & i>j.\sluttuborg
\end{eqnarray}
This structure makes $\Delta_X ^{\bullet} = X \times_k \Delta^{\bullet}$ a cosimplicial scheme.

\begin{defi}{\rm Let $\calZ_q (X, n)$ be the free abelian group on irreducible closed subvarieties in $\Delta_X ^n$ of dimension $q$ that intersect with all faces properly, i.e. in the right codimensions.
}
\end{defi}

The faces and degeneracies induce maps
\begin{eqnarray}
\partial_i : &\calZ_q (X, n) \to \calZ_{q-1} (X, n-1), \ \ 0 \leq i \leq n \\
\pi_j :& \calZ_q (X, n) \to \calZ_{q+1} (X, n+1), \ \ 0 \leq j \leq n
\end{eqnarray} making the groups $\calZ_q (X, n)$, a simplicial abelian group. The associated complex \begin{equation}
\cdots \overset{\partial}{\to} \calZ_{q+1} (X, n+1) \overset{\partial}{\to} \calZ_q (X, n) \overset{\partial}{\to} \calZ_{q-1} (X, n-1) \overset{\partial}{\to} \cdots
\end{equation} is the higher Chow complex. For $n \geq 0$, the \emph{higher Chow groups} are defined to be the homology groups of the above complex. They are written as $CH_q (X, n)$. In case $X$ is equidimensional, we also write
\begin{equation}
CH^p (X, n):= CH_{n + \dim X - p } (X, n).
\end{equation}

For the polynomial ring $k[t_0, \cdots, t_n]$ in $(n+1)$-variables, define 
\begin{eqnarray}
&T^{-1}=T_n ^{-1}: k[t_0, \cdots, t_n ] \to k[t_0, \cdots, t_n] \\
& t_i \mapsto \tuborg t_{i+1} & 0 \leq i \leq n-1,\\ t_0 & i=n. \sluttuborg
\end{eqnarray} It descends to $
T_n  ^{-1}: k[t_0, \cdots, t_n] / \left(\sum  t_i -1\right) \to k[t_0, \cdots, t_n] /\left( \sum t_i -1 \right)$ thus induces a map $T_n$ on the $k$-rational points of $\Delta^n$ sending $(t_0, \cdots, t_n) \mapsto (t_n, t_0, \cdots, t_{n-1})$ and further determines a map
\begin{equation}
T_n: \calZ_q (X, n) \to \calZ_q (X, n).
\end{equation}

\begin{prop}
It satisfies the cyclic compatibility conditions \eqref{cyclic compatibility}.
\end{prop}

\begin{proof}After a moment's thought, we can see that it is enough to prove it for $X= \spec (k)$ and for $k$-rational points of $\Delta_n$. For this case, this is apparent. But, we can carry out a different rigorous proof on the level of the simplicial ring $\frac{k[t_0, \cdots, t_{n}]}{ \left(\sum t_i -1\right)}$ with the face maps
\begin{equation}
\partial_j ^* (t_i) = \tuborg t_i & i<j, \\ 0  & i=j, \\ t_{i-1} & i>j. \sluttuborg
\end{equation} The identities we have to prove are
\begin{equation}
\tuborg \partial_0 ^* T^{-1} = \partial_n ^* & \\ \partial_j ^* T^{-1} = T^{-1} \partial_{j-1} ^* & 1 \leq j \leq n.\sluttuborg
\end{equation}

But, by straightforward calculations,
\begin{eqnarray}
\partial_0 ^* T^{-1} (t_i) &=& \tuborg \partial_0 ^* t_{i+1} & 0 \leq i \leq n-1 \\ \partial_0 ^* t_0 & i=n \sluttuborg \\
&=& \tuborg t_i & 0 \leq  i \leq n-1 \\ 0 & i=n \sluttuborg \\
& = & \partial_n ^* (t_i)
\end{eqnarray} thus $\partial_0 ^* T^{-1} = \partial_n ^*$, and for $1 \leq j \leq n$,

\begin{eqnarray}
\partial_j ^* T^{-1} (t_i) &=& \tuborg \partial_j ^* t_{i+1} & 0 \leq i \leq n-1 \\ \partial_j ^* t_0 & i=n \sluttuborg \\
&=& \tuborg  t_{i+1} & i+1 <j \\ 0 & i+1 = j \\ t_i & i+1 >j \\ t_0 & i=n \sluttuborg \\
T^{-1} \partial_{j-1} ^* (t_i) &=& \tuborg T^{-1} t_i & i<j-1 \\ 0 & i=j-1 \\ T^{-1} t_{i-1} & i>j-1 \sluttuborg \\
&=& \tuborg  t_{i+1} & i+1 <j \\ 0 & i+1 = j \\ t_i & i+1 >j \\ t_0 & i=n \sluttuborg
\end{eqnarray} which shows that $\partial_j ^* T^{-1} = T^{-1} \partial_{j-1} ^*$. This finishes the proof.\end{proof}

The cyclic homology of the cyclic object $\calZ_q (X, n)$ will be denoted by $CHC_q (X, n)$ and it will be called the \emph{Cyclic higher Chow group}. Certainly these groups form the Connes periodicity sequence together with the higher Chow groups. When $X$ is equidimensional, we define $CHC^p (X, n) = CHC_{n + \dim X - p}(X, n)$. A care should be taken here, however. Even though we use this (co)-dimension notation, the group $CHC_q (X, n)$ may contain classes represented by cycles of dimension $\leq q$ by its definition. Thus, the natural map $CHC_q (X, n) \to CHC_q ^{\lambda} (X, n)$, to the \emph{Connes higher Chow group} cannot be an isomorphism, in general. But, this is not too far from being an isomorphism as we have seen previously. Namely, the cycles in the Connes higher Chow groups of dimension $<q$ are all torsion and up to $\otimes \bbQ$, these groups are isomorphic. We will record it here as an independent theorem, even though it doesn't require a new proof.

\begin{thm}\label{Cyclic to Connes}The natural map induces an isomorphism
\begin{equation}CHC_q (X, n) \otimes_{\bbZ} \bbQ \overset{\simeq}{\to} CHC^{\lambda} _q (X, n) \otimes_{\bbZ} \bbQ.\end{equation}
\end{thm}

\begin{remark2}{\rm The cyclic operator $t$ does not act nicely on the (Cyclic) higher Chow groups in general. However, by definition, on the Connes higher Chow groups it acts as $(-1)^n$ thus the Connes higher Chow groups are nice $\bbZ/(n+1)$-modules. To see why Connes higher Chow groups can be useful in studies of motives, see Theorem \ref{equivalence}.
}
\end{remark2}

\begin{remark2}{\rm Proposition (1.3) in \cite{B1} shows that the complex $(\calZ_{\bullet} (X, \cdot), \partial)$ is covariant functorial for proper maps and contravariant funtorial for flat maps, thus the higher Chow groups have this functorial property. Using the same argument, we can also prove that the acyclic complex $(\calZ_{\bullet}(X, \cdot), \partial')$ shares the same functorial property. Hence so does the cyclic bicomplex $CC(X,r)$ and the Cyclic higher Chow groups have the same functorial property. If $X \to Y$ is any morphism between smooth varieties, then the Cyclic higher Chow groups are functorial.}\end{remark2}

More properties of the Cyclic higher Chow groups will be investigated in a later section.

\subsection{Additive higher Chow complex}

For a number $r \not = 0$, 
\begin{eqnarray*} & & \Delta^n  =  \spec \left( k[t_0, \cdots, t_n]/( \sum t_i - 1 ) \right) \\ & \simeq & \Delta_r ^n  = \spec \left( k[t_0, \cdots, t_n] / (\sum t_i -r ) \right).\end{eqnarray*} By letting $r \to 0$, more precisely, by considering $$Q^n = \spec \left( k[t_0, \cdots, t_n]/ (\sum t_i ) \right),$$ we can obtain an additive version of the previous higher Chow theory. For the sake of convenience of readers, let us recall some necessary definitions from \cite{BE1}. Let $k$ be a field and $X$ be a scheme of finite type over $k$. Let $$Q^n = \spec \left( k[t_0, \cdots, t_n]/ \sum_{i=0} ^n t_i \right)$$ with the faces

\begin{eqnarray}
\partial_j: Q^{n-1} \to Q^n ; \ \ \ \partial_j ^* (t_i) = \tuborg t_i & i<j, \\ 0 & i=j,\\ t_{i-1} & i>j.\sluttuborg
\end{eqnarray}
and the degeneracies
\begin{eqnarray}
\pi_j: Q^n \to Q^{n-1} ;\ \ \  \pi_j ^* (t_i) = \tuborg t_i & i<j, \\ t_i + t_{i+1} & i=j, \\ t_{i+1} & i>j.\sluttuborg
\end{eqnarray}
This structure makes $Q_X ^{\bullet} = X \times_k Q^{\bullet}$ a cosimplicial scheme.
\begin{defi}{\rm Let $S\calZ_q (X, n)$ be the free abelian group on irreducible closed subvarieties in $Q_X ^n$ of dimension $q$ with the property:
\begin{enumerate}
\item [(i)] they don't meet $X \times \{ 0 \}$, and
\item [(ii)] they meet all faces properly.
\end{enumerate}
}
\end{defi}

The faces induce maps
\begin{equation}
\partial_i : S\calZ_q (X, n) \to S\calZ_{q-1} (X, n-1), \ \ 0 \leq i \leq n
\end{equation} thus with $\partial = \sum_{i=0} ^n (-1)^i \partial_i$, $S\calZ_{q + \bullet} (X, \bullet)$ is a complex of abelian groups
\begin{equation}
\cdots \overset{\partial}{\to} S\calZ_{q+1} (X, n+1) \overset{\partial}{\to} S\calZ_q (X, n) \overset{\partial}{\to} S\calZ_{q-1} (X, n-1) \overset{\partial}{\to} \cdots
\end{equation}
We will call it the additive higher Chow complex. For $n \geq 1$, the \emph{additive higher Chow groups} are defined to be the homology groups of the above complex. They will be written as $SH_q (X, n)$. In case $X$ is equidimensional, we will also write
\begin{equation}
SH^p (X, n):= SH_{n + \dim X - p } (X, n).
\end{equation}

The rest of the discussion works in exactly same way for the same operator $T^{-1}$. We will not repeat it here again. Some properties of the additive Cyclic higher Chow groups will be mentioned in the next section.

\section{Properties of Cyclic higher Chow groups}

\subsection{Properties}Various properties of the Cyclic higher Chow groups can be easily transferred from either the properties of the higher Chow groups using the Connes periodicity sequence, or by adapting the methods used to prove them. Properties of the higher Chow groups recalled here can be found in \cite{B1} and \cite{BM}.

\begin{cor}Let $k$ and $X$ be as above. Then,
\begin{eqnarray} & CHC_0 (X, n) \simeq CH_0 (X, n) \\
& CHC_q (X, 0) \simeq CH_q (X, 0) \simeq CH_q (X).
\end{eqnarray}
\end{cor}
\begin{proof}The tails of the Connes periodicity sequences give the isomorphisms:
\begin{equation}
\cdots \to CHC_1 (X, n+1) \to 0 \to CH_0 (X, n) \overset{I}{\to} CHC_0 (X, n) \to 0,
\end{equation}
\begin{equation}
\cdots \to CHC_{q+1} (X, 1) \to 0 \to CH_q (X, 0) \overset{I}{\to} CHC_q (X, 0) \to 0.
\end{equation}
\end{proof}

\begin{prop}[codimension 1]Let $X$ be a smooth variety over $k$. Then,
\begin{equation}
CHC^1 (X, n) = \tuborg CH^1 (X, 0) = \pic (X) & n :\mbox{ even}, \\
CH^1 (X, 1) = \Gamma(X, \calO_X ^{\times}) & n :\mbox{ odd}. \sluttuborg
\end{equation}
\end{prop}
\begin{proof}
Notice that for a fixed Connes periodicity exact sequence, the codimensions of all involved groups are same. In particular, we can look at the Connes sequence of codimension $1$. Recall (Theorem (6.1) in \cite{B1}) that
\begin{equation}
CH^1 (X, n) = \tuborg \pic (X) & n=0, \\ \Gamma(X, \calO_X ^{\times}) & n=1, \\ 0 & n \geq 2.\sluttuborg
\end{equation} For $n=0$, we already know $CHC^1 (X,0)$. For $n\geq 1$, the part of the Connes sequence
\begin{equation}
\cdots \to CH^1 (X, n+3) \overset{I}{\to} CHC^1 (X, n+2) \overset{S}{\to} CHC^1 (X, n) \overset{B}{\to} CH^1 (X, n+1)\end{equation} shows that $CHC^1 (X, n+2) \simeq CHC^1 (X, n)$. Near the tail of the sequence, we obtain
\begin{equation}
0  \overset{I}{\to} CHC^1 (X, 2) \overset{S}{\to} CHC^1 (X, 0) \overset{B}{\to} CH^1 (X, 1) \overset{I}{\to} CHC^1 (X, 1) \to 0.
\end{equation} The map $B: CHC^1 (X,0) \to CH^1 (X, 1)$ is $H^1 (X, \calO_X ^{\times}) \to H^0 (X, \calO_X ^{\times})$ is a connecting homomorphism and it is $0$ by a simple application of the serpent lemma. Thus, $CHC^1 (X, 2) \simeq CHC^1 (X, 0) \simeq \pic (X)$ and $\Gamma(X, \calO_X ^{\times}) \simeq CH^1 (X, 1) \simeq CHC^1 (X, 1)$.
\end{proof}

\begin{thm}\label{equivalence}Let $f: X \to Y$ be a good morphism for which the higher Chow groups and the Cyclic higher Chow groups are functorial. 

Then, it induces isomorphisms of the higher Chow groups if and only if it induces isomorphisms of the Cyclic higher Chow groups. 

In particular, $f$ induces isomorphisms of the rational motivic cohomology groups if and only if it induces isomorphisms of the rational Connes higher Chow groups.
\end{thm}

\begin{proof}For $0$-cycles, the above corollary answers the question. We then use an induction argument (\emph{c.f.} Corollary 2.2.3 in \cite{L}) on the dimension of the cycles. Suppose that the morphism induces isomorphisms on the higher Chow groups. Then, the $5$-lemma on the part of the periodicity sequence
\begin{equation}
\cdots \to CHC_0 (\cdot, *-1) \to CH_1 (\cdot, *) \to CHC_1 (\cdot, *) \to 0 \to CH_0 (\cdot, *-1) \to \cdots
\end{equation} shows that the morphism induces isomorphisms on $CHC_1 (\cdot, *)$. Suppose that the morphism induces isomorphisms for $CHC_i (\cdot, *)$ for $i \leq r-1$. Then, by applying the $5$-lemma again on 
\begin{eqnarray*}
&\cdots \to CHC_{r-1} (\cdot, *-1) \overset{B}{\to} CH_r (\cdot, *)\overset{I}{ \to} CHC_r (\cdot, *) \overset{S}{\to} \\ &CHC_{r-2} (\cdot, *-2) \overset{B}{\to} CH_{r-1} (\cdot, *-1) \to \cdots
\end{eqnarray*} we see that the morphism induces isomorphisms for $CHC_r (\cdot, *)$. Thus by induction we have the desired result. The proof for the other direction can be done similarly. The last remark concerning rational motivic cohomology groups and the Connes higher Chow groups follows from \cite{Be2}, Theorem (9.1) in \cite{B1} and the Theorem \ref{Cyclic to Connes}.
\end{proof}

\begin{cor}[$\bbA^1$-homotopy]The Cyclic higher Chow groups are $\bbA^1$-homotopy invariant. More precisely, if $E$ is a vector bundle on $X$, then $CHC^* (X, n) \simeq CHC^* (E, n)$.
\end{cor}
\begin{proof} Because the higher Chow groups are $\bbA^1$-homotopy invariant (See Theorem (2.1) in \cite{B1}.), by the Corollary \ref{equivalence}, we must have the same property for the Cyclic higher Chow groups.
\end{proof}

\bigskip Let $X$ be a quasiprojective algebraic $k$-scheme and let $Y \subset X$ be a closed subscheme. Let $U= X-Y$. From \cite{BM}, we know the moving lemma that says: 
\begin{equation}\label{moving lemma}
\calZ_{\bullet} (X, *)/ \calZ_{\bullet}(Y, *) \to \calZ_{\bullet}(U, *)
\end{equation}
is a homotopy equivalence. This immediately gives the localization long exact sequence for the higher Chow groups:
\begin{eqnarray*}
&\cdots \to CH_{q+1}(U, n+1) \to CH_q(Y, n) \to CH_q (X, n) \to \\ & CH_q(U, n) \to CH_{q-1}(Y, n-1)\to \cdots
\end{eqnarray*}

Similar results can be proven for the Cyclic higher Chow groups and the Connes higher Chow groups. 

\begin{prop}[Localization for Connes groups] Under the above assumptions on $Y, X$ and $U$, there is a localization long exact sequence for the Connes higher Chow groups
\begin{eqnarray*}
& \cdots \to CHC^{\lambda}_{q+1}(U, n+1) \to CHC^{\lambda}_q(Y, n) \to CHC^{\lambda}_q (X, n) \to \\ & CHC^{\lambda}_q(U, n) \to CHC^{\lambda}_{q-1}(Y, n-1)\to \cdots.
\end{eqnarray*}
\end{prop}
\begin{proof}The sequence
\begin{equation}
\calZ_{\bullet} (Y, *) \to \calZ_{\bullet}(X, *) \to \calZ_{\bullet}(U, *)
\end{equation} is $(1-t)$-equivariant so that the induced sequence of the Connes higher Chow complexes
\begin{equation}
\calZ_{\bullet}(Y, *)/(1-t) \to \calZ_{\bullet}(X, *)/(1-t) \to \calZ_{\bullet}(U, *)/ (1-t)
\end{equation} gives a quasi-isomorphism
\begin{equation}
\frac{ \calZ_{\bullet} (X,*)/ (1-t)}{\calZ_{\bullet}(Y,*)/(1-t)} \to \calZ_{\bullet}(U, *)/(1-t). 
\end{equation} This proves that the localization sequence is exact.\end{proof}

From the Theorem \ref{Cyclic to Connes}, up to $\otimes \bbQ$ the Connes higher Chow groups are isomorphic to the Cyclic higher Chow groups. Thus, the rational Cyclic higher Chow groups also have the localization exact sequence. But actually it can be proven even without taking $\otimes \bbQ$.

\begin{prop}[Localization for Cyclic groups]Under the above assumptions on $Y, X$ and $U$, there is a localization long exact sequence for the Cyclic higher Chow groups
\begin{eqnarray*}
 &\cdots \to CHC_{q+1}(U, n+1) \to CHC_q(Y, n) \to \\ &CHC_q (X, n) \to CHC_q(U, n) \to CHC_{q-1}(Y, n-1)\to \cdots.
\end{eqnarray*} Furthermore, this sequence is compatible with the Connes periodicity sequence in the following sense: we have a commutative diagram with exact rows, that are localization sequences, and exact columns, that are the Connes periodicity sequences:
$$
\xymatrix{  & \vdots \ar[d]_B & \vdots \ar[d]_B & \vdots \ar[d] _B & \\
 \ar[r] & CH_q (Y,n) \ar[r] \ar[d]_I & CH_q (X, n) \ar[r] \ar[d]_I & CH_q (U, n) \ar[r] \ar[d]_I&  \\
 \ar[r] &CHC_q (Y, n) \ar[r] \ar[d]_S & CHC_q (X, n) \ar[r] \ar[d] _S & CHC_q (U, n) \ar[r] \ar[d]_S  &  \\
 \ar[r] &CHC_{q-2} (Y, n-2) \ar[r] \ar[d]_{B} & CHC_{q-2} (X, n-2) \ar[r] \ar[d]_{B} & CHC_{q-2} (U, n-2) \ar[r] \ar[d]_{B}&  \\
 \ar[r] & CH_{q-1} (Y, n-1) \ar[r] \ar[d]_I & CH_{q-1} (X, n-1) \ar[r] \ar[d]_I & CH_{q-1} (U, n-1) \ar[r] \ar[d]_I &  \\
& \vdots & \vdots & \vdots  & }
$$
\end{prop}

\begin{proof}The proof is also easy. The homotopy equivalence $$\calZ_{\bullet}(X, *) / \calZ_{\bullet}(Y, *) \to \calZ_{\bullet} (U, *)$$ induces a quasi-isomorphism
\begin{equation}
CC(X, r)/ CC(Y, r)\to CC(U, r)
\end{equation} which gives the localization sequence for the Cyclic higher Chow groups. The compatibility can be seen by looking at the commutative diagram of bicomplexes with exact rows and columns and applying the moving lemma to the bottom row:
$$
\xymatrix{ & 0 \ar[d] & 0 \ar[d] & 0 \ar[d] & \\
0 \ar[r] & CC(Y,r)^{\{2 \}} \ar[r] \ar[d]  & CC(Y, r) \ar[r] \ar[d] & CC(Y, r)[2,0] \ar[r] \ar[d] & 0 \\
0 \ar[r] & CC(X, r)^{\{2 \}} \ar[r] \ar[d] & CC(X, r) \ar[r] \ar[d] & CC(X, r)[2,0] \ar[r] \ar[d] & 0 \\
0 \ar[r] & CC(U, r) ^{\{2 \}} \ar[r] & CC(U, r) \ar[r] & CC(U, r) [2, 0] \ar[r] & 0 }
$$
\end{proof}

\begin{prop}[Local to global spectral sequence]Let $\calT ot_{\bullet}$ be the complex of Zariski sheaves on $X$ concentrated in negative degrees associated to the presheaves $U \mapsto Tot_{\bullet}(CC(U, r))$ for Zariski open $U \subset X$. Then, $CHC_* (X, n) = \bbH^{-n} (X, \calT ot_{\bullet})$. In particular, there is a spectral sequence 
\begin{equation}\label{spectral sequence}
E_2 ^{p,q} = \bbH^p (X, \calC\calH\calC_s (-q) ) \Rightarrow CHC_s (X, -p -q)
\end{equation} where $\calC\calH\calC_s (n)$ is the Zariski sheaf associated to the presheaf $$U \mapsto CHC_s (U, n).$$ Furthermore, this spectral sequence is compatible with the Conne periodicity.
\end{prop}
\begin{proof} From the Theorem (3.4) in \cite{B1}, we know that the presheaf $\calS_n$ on $X$ defined by $U \mapsto \calZ_* (X, n) / \calZ_* (X-U, n)$ is in fact a flasque sheaf. Thus, the presheaf $\calT_n$ on $X$ defined by $$U \mapsto Tot_n (CC(X,r))/ Tot(CC(X-U, r))$$ is also a flasque sheaf. By the moving lemma \eqref{moving lemma}, the natural morphism of complexes of sheaves $\calT_{\bullet} \to \calT ot_{\bullet}$ is a quasi-isomorphism and since each $\calT_n$ is flasque, we have $\bbH^{-n} (X, \calT ot_{\bullet}) = \bbH^{-n} (X, \calT_{\bullet}) = H_n (Tot_n (CC(X, r))) = CHC_* (X, n)$. Since $\calC \calH \calC_* = H^{-n} (\calT_{\bullet})$, the standard spectral sequence gives \eqref{spectral sequence}. Compatibility with the Connes sequence is an easy exercise.\end{proof}

\begin{remark2}{\rm Some results on the higher Chow groups related to the D. Quillen's higher algebraic $K$-theory (See \cite{Q}), like analogues of Chern classes, finite coefficients, Gersten's conjecture and projective bundle theorem, could not be transplanted on this ground yet as we lack a cyclic analogue of the Quillen's higher $K$-theory. Will there be an easy way to define such an object?
}
\end{remark2}

\subsection{Conjectures}It is possible to translate several conjectures on the motivic cohomology groups in terms of their cyclic cousins. Let $X$ be a smooth quasi-projective variety over a field $k$. Recall that the integral motivic cohomology groups are isomorphic to the higher Chow groups and the relations (See \cite{V}) are
\begin{equation}
CH^p(X, n) \simeq H_{\calM} ^{2p-n} (X, \bbZ(p)).
\end{equation}

One famous hard conjecture on these groups is the following:

\begin{conj}[Beilinson-Soul\'e vanishing (See \cite{VSF})] Under the above assumption, for all $i<0$ \begin{equation}\label{BS vanishing}H_{\calM} ^i (X, \bbZ(p)) = 0.\end{equation}
\end{conj}

In terms of higher Chow groups, it is equivalent to the statment that $CH^p (X, n)= 0$ for all $n > 2p$. The Connes periodicity sequence supplies an equivalent formulation of the above conjecture in terms of the Cyclic higher Chow groups:

\begin{thm}\label{reformulation}The Beilinson-Soul\'e vanishing conjecture is true if and only if the following conditions hold:
\begin{enumerate}
\item via the shift map $S$, $CHC^p (X, n+1) \overset{\simeq}{\to} CHC^p (X, n-1)$ for $n > 2p$.
\item the shift map $CHC^p (X, 2p+1) \overset{S}{\to} CHC^p (X, 2p-1)$ is injective.
\end{enumerate}
\end{thm}

\begin{proof}This is a simple consequence of the Connes periodicity sequence
\begin{eqnarray*}
& \cdots \to CHC^p (X, n-1) \overset{B}{\to} CH^p (X, n) \overset{I}{\to}\\ & CHC^p (X, n) \overset{S}{\to} CHC^p (X, n-2) \overset{B}{\to} CH^p (X, n-1)\to \cdots
\end{eqnarray*} and to note \eqref{BS vanishing}.\end{proof}

\begin{remark2}{\rm Up to $\otimes \bbQ$, by the Theorem \ref{Cyclic to Connes}, the Connes higher Chow groups are isomorphic to the Cyclic higher Chow groups. We know that the Connes higher Chow group $CHC_q ^{\lambda} (X, n)$ has a natural action of $\bbZ/ (n+1)$ and its generator acts as the sign change $(-1)^n$. Thus, the first part of the Theorem \ref{reformulation} implies that when $n >2p$, there are certain abelian groups equipped with the both actions of $\bbZ/n$ and $\bbZ/(n+2)$ and both of them act as sign changes.}
\end{remark2}

Now, suppose that $X$ is a scheme of finite type over $\bbZ$. In \cite{B1}, under the assumptions that
\begin{enumerate}
\item $CH^r (X, i)$ can be defined for schemes $X$ of finite type over $\bbZ$,
\item they have finite ranks, and
\item $CH^r (X, i)\otimes \bbQ = 0$ for $i >>0$,
\end{enumerate}
S. Bloch reformulated the conjecture (See \cite{So}) of C. Soul\'e on the order of vanishing of the zeta function $\zeta_X (s)$ as follows:

\begin{conj}[Soul\'e]Let $d = \dim X$. Then,
\begin{equation}
- \ord_{s= d -r } \zeta_X (s) = \sum_i (-1)^i rk (CH^r (X, i)).
\end{equation}
\end{conj}

We can also reformulate this conjecture in terms of the Cyclic higher Chow groups. If $CH^r(X, i)$ can be defined for the above schemes $X$, certainly we should be able to define $CHC^r(X,i)$. Furthermore, by the $5$-lemma applied to the Connes periodicity sequence, we see that $CH^r (X, i)$ have finite ranks if and only if $CHC^r (X,i)$ have finite ranks. Thus, we have

\begin{thm}Suppose that there is an integer $N$ such that  $rk(CH^r (X, n)) = 0$ for all $n \geq N$. (Notice that this assumption is a weaker version of the the Beilinson-Soul\'e vanishing conjecture for $CH^r (X, n)$.) Then, the conjecture of Soul\'e is true if and only if
\begin{equation}
\ord_{s= d-r} \zeta_X (s) =(-1)^{N} rk(CHC ^r (X, N-1)) +(-1)^{N+1} rk (CHC^r (X, N)).
\end{equation}
\end{thm}

\begin{remark2}{\rm The quantity on the right hand side appears to be depending on $N$, but in fact it doesn't: the Connes periodicity sequence shows that 
\begin{eqnarray*}
\cdots \to CH^r (X, N+1) \to CHC^r (X, N+1) \to \\ CHC^r (X, N-1) \to CH^r (X, N) \to \cdots
\end{eqnarray*}is exact and the first and the fourth entry having zero ranks, two middle groups have the same ranks. Thus, 
\begin{eqnarray*}
&(-1)^{N} rk(CHC ^r (X, N-1)) +(-1)^{N+1} rk (CHC_r (X, N)) \\
&= (-1)^{N+1} rk(CHC_r (X, N)) + (-1)^{N+2} rk (CHC^r(X, N+1))
\end{eqnarray*} and hence, by induction we have the equality for any number $n \geq N$. If an analogue of the Beilinson-Soul\'e vanishing conjecture is true for these schemes, then $CH^r (X, n) =0$ for $n > 2r$. Thus certainly $N= 2r+1$ should work and the above conjecture is equivalent to
\begin{equation}\ord_{s=d-r} \zeta_X (s) = - rk (CHC^r (X, 2r)) + rk (CHC^r (X, 2r+1)),\end{equation}which is simpler than the original conjecture of C. Soul\'e.
}
\end{remark2}

\begin{proof}The proof is again a simple consequence of the Connes periodicity sequence. Since $rk(CH^r (X, N))=0$, after tensoring with $\bbQ$, we can truncate the long exact sequence at this point and look at only the right hand side of this sequence. In particular, we can notice that for $i < N-1$, $CHC^r (X, i)$ appears twice and their degrees differ by $5$ in the long exact sequence. Thus when we take the alternating sum of ranks, ranks of the groups $CHC^r (X, i)$ with $i<N-1$ will be cancelled and only the ranks of $CHC^r (X, N-1)$, $CHC^r (X, N)$ and $CH^r (X, i)$, $i<N$ are not cancelled. Hence, we have the equation
\begin{eqnarray*}
\sum_{i=0} ^{N-1} rk (CH^r (X, i)) + (-1)^N rk (CHC^r (X, N-1)) \\ + (-1)^{N+1} rk (CHC^r (X, N)) = 0.
\end{eqnarray*}This proves the theorem.
\end{proof}

\subsection{The additive case} As applications of the Connes periodicity sequence we can readily deduce some facts about the additive Cyclic higher Chow groups from facts about the additive higher Chow groups.

\begin{cor}\begin{eqnarray}
&SHC_0 (X, n) \simeq SH_0 (X, n), n \geq 1 \\
& SHC_q (X, 1) \simeq SH_q (X, 1)
\end{eqnarray}
\end{cor}
\begin{proof}Easy exericise.
\end{proof}

\begin{thm}Let $f: X \to Y$ be a morphism that induces good functorial properties on $SH$ and $SHC$. Then $f$ induces isomorphisms on the additive higher Chow groups if and only if it induces isomorphisms on the additive Cyclic higher Chow groups.
\end{thm}
\begin{proof}Easy exercise.
\end{proof}

\begin{prop}[codimension 1] \begin{equation}
SHC ^1 (k, n) \simeq \tuborg SH^1 (k, 1) \simeq k & n = \mbox{odd}, \\ 0 & n = \mbox{even}.\sluttuborg
\end{equation}
\end{prop}
\begin{proof}Recall (Proposition 7.2 in \cite{BE1}) that \begin{equation} SH^1 (k, n) = \tuborg k & n=1, \\ 0 & n \geq 2.\sluttuborg\end{equation}

We already know that $SHC^1 (k,1) \simeq SH^1 (k,1)$. For $n=2$, the part of the Connes periodicity sequence
\begin{equation}
0 = SH^1 (k,2) \overset{I}{\to} SHC^1 (k,2) \overset{S}{\to} 0
\end{equation} shows that $SHC^1 (k,2) \simeq 0$. In general, for $n \geq 1$, the Connes exact sequence
\begin{equation}
 0=SH^1 (k, n+2) \overset{I}{\to} SHC^1 (k, n+2) \overset{S}{\to} SHC^1 (k, n) \overset{B}{\to} SH^1 (k, n+1)=0
\end{equation}shows that $SHC^1 (k, n+2) \simeq SHC^1 (k, n)$. This finishes the proof.\end{proof}

\bigskip

\noindent \textbf{Acknowledgement} The author would like to thank Alexander Beilinson, Spencer Bloch, Andrew Blumberg, Peter May, Madhav Nori, and Shmuel Weinberger for their comments on this work. This paper is based on a chapter in the author's doctoral thesis at the University of Chicago.

\end{document}